\documentclass[preprint,12pt]{elsarticle}

\usepackage{graphicx}
\usepackage{amssymb}

\usepackage{epsf}
\usepackage{color}
\usepackage[dvips]{epsfig}
\usepackage{natbib}
\usepackage{amsbsy}
\usepackage{amsfonts}
\usepackage{amsmath}
\usepackage{latexsym}
\usepackage{rotating}
\usepackage{geometry}
\usepackage{subfigure}

\usepackage{array}
\usepackage{multirow}
\usepackage{tabularx,tabulary}
\usepackage{rotating}

\usepackage{pst-all}
\usepackage{pstricks}
\usepackage{pstricks-add}
\usepackage{pst-lens}
\usepackage{pst-plot}
\usepackage{acronym}

\newcommand{\ga}{{jetG1}}

\newcommand{\gb}{{jetG2}}
\newcommand{\gc}{{jetG3}}
\newcommand{\gd}{{jetG4}}

\newcommand{\jetLES}{{LES3D}}
\newcommand{\jetELES}{{ELES3D}}
\newcommand{\jetDNS}{{DNS3D}}
\newcommand{\PD}[2]{{ {\partial#1}\over {\partial #2}}}
\newcommand{\jettwo}{{DNS2D}}

\newcommand{\gapCtrl}[2]{gap#1F#2}

\journal{arXiv}

\begin{document}

\begin{frontmatter}


\title{On the Stability of Gradient Based Turbulent Flow Control without Regularization}



\author{Emre \"Ozkaya, Nicolas R. Gauger}

\address{Chair of Scientific Computing, TU Kaiserslautern, Paul-Ehrlich-Strasse 34, 67663 Kaiserslautern, Germany}
\author{Daniel Marinc\footnote{Former Ph.D. student at the Chair of Fluid Mechanics, University of Siegen, 57062 Siegen, Germany}, Holger Foysi}

\address{Chair of Fluid Mechanics, University of Siegen, Paul-Bonatz-Str. 9-11, 57062 Siegen, Germany}

\begin{abstract}
In this paper, we discuss selected adjoint approaches for the turbulent flow control. In particular, we focus on the application of adjoint solvers for the scope of noise reduction, in which flow solutions are obtained by large eddy and direct numerical simulations. Optimization results obtained with round and plane jet configurations are presented. The results indicate that using large control horizons poses a serious problem for the control of turbulent flows due to existence of very large sensitivity values with respect to control parameters. Typically these sensitivities grow in time and lead to arithmetic overflow in the computations. This phenomena is illustrated by a sensitivity study performed with an exact tangent-linear solver obtained by algorithmic differentiation techniques. 
\end{abstract}


\end{frontmatter}


\section{Introduction}
\label{Sec:Intro}

The control of turbulent flows has been an area
of particular interest in aerospace research due to major commercial benefits. For most applications, the control parameters are determined by using simulation based optimization techniques. As typical examples, active flow control of high-lift devices for lift enhancement \citep{nielsen_2011}, flow relaminarization \citep{Bewley:2001} and noise control \citep{Wei:2006} can be given. 

Especially for large-scale problems, an efficient way of finding the optimal values of control parameters is by employing gradient-based numerical optimization techniques. These methods are ideally combined with continuous or discrete adjoint solvers to evaluate the gradient. The main advantage of the adjoint approaches is that they enable efficient evaluation of gradients, irrespective of the number of control parameters, at a fixed computational cost. Therefore, numerical optimization studies using high dimensional control vectors and high fidelity simulation for function evaluations become viable within limited computational resources. 

Broadly, we can classify adjoint approaches
into two categories: continuous and discrete adjoint methods. In the continuous adjoint method \citep{Jameson:1998}, 
the optimality system is derived from the continuous optimization problem. The resulting adjoint partial differential equations (PDEs) are then discretised and solved using state-of-the-art numerical
methods. Although being computationally efficient, development of continuous adjoint flow solvers requires considerable development effort and a good understanding of the underlying flow solver is often required. Furthermore, their maintenance becomes problematic as the underlying non-linear flow solvers are subject to continuous modifications, e.g., new boundary conditions, new physical models etc. In short, extension of the adjoint solver to incorporate the new features might be a quite challenging task. 

In the discrete adjoint method, on the other hand, the derivation of the optimality conditions starts directly with the discretized state PDEs that govern the fluid flow. Based on a given discretization, the discrete adjoint equation is derived. In general, compared to the continuous adjoint solvers, discrete adjoint solvers are more straightforward to implement. Therefore, they have found a wider acceptance for flow control applications of practical relevance in the past.

In general, a discrete adjoint method for optimal active flow control can be developed either by
using the so-called hand-discrete approach \citep{Nielsen} or by employing Algorithmic Differentiation (AD) techniques to the underlying non-linear flow solver. In the hand-discrete approach, the adjoint equations are derived by linearizing the discrete residuals by hand. Based on the derivation, a computer code is then implemented to solve the adjoint equations and to evaluate the gradient vector. In the AD based approach, on the other hand, the adjoint code is generated directly by applying AD techniques to the computational fluid dynamics (CFD) code that is implemented to solve the discretized flow equations. 

Accurate computation of sensitivities requires exact differentiation of all terms that constitute the discrete residual. However, exact linearization of these terms might be often quite complex, laborious and error prone. To simplify this tedious effort and ease the development of the adjoint solvers, various Jacobian approximations have been proposed in the past \citep{Pet2010}. While using these approximations, linearization of certain terms in the flux Jacobian is omitted. As typical examples, turbulent models, flux limiters, convergence acceleration schemes and higher order reconstruction terms can be given. The Jacobian approximations ease the development effort significantly at the expense of some inaccuracy in the sensitivity evaluation. Usually for steady-state problems, the inaccuracies incurred by the Jacobian approximations can be tolerated. In unsteady flows, however, the effect of these approximations on the accuracy of sensitivities is much more significant as the errors generated in the adjoint solution tend to accumulate rapidly while solving the adjoint equations backward-in-time \citep{anil2013}. On the other hand, by using the AD techniques in reverse mode, very accurate adjoint solvers can be generated since the exact differentiation of all residual terms can be performed by AD packages with much ease. Therefore, Jacobian approximations are no more required in the AD based discrete adjoint framework.

The AD based adjoint solvers, apart from the advantages mentioned above, have in general higher memory demand compared to the continuous and hand-discrete counterparts. The reason for that is the fact that all the operations that are performed within a time iteration must be saved while integrating the state vector forward in time. Especially for large eddy (LES) or direct numerical (DNS) simulations, the number of floating point operations can be very large, therefore the memory demand might required for the reversal of a time iteration might overwhelmingly large. In addition, extra memory is also required for the reversal of the time loop. In the most extreme case, the complete forward trajectory of the state solution is kept in memory. This approach is known as the store-all approach \citep{HascoetGK09} in the AD community. As far as the run-time requirements are concerned, for most applications, the store-all approach is not feasible. A simple way of overcoming the excessive memory demand of the store-all approach is by storing parts or the entire forward trajectory on the hard-disk at the expense of increased run-time. This approach has been widely used in the past \citep{Mani08,Nada02, Rump07} for unsteady adjoint computations. Another memory saving solution is the checkpointing strategy. In adjoint solvers based on checkpointing, the flow solutions are stored only at selective time iterations known as checkpoints. Various checkpointing strategies have been proposed based on the storage criteria for checkpoints \citep{Griewank2000,Stumm09}. 
Although, checkpointing strategies are successfully applied to for unsteady Reynold-averaged Navier-Stokes (URANS) and LES computations recently \citep{Zhou_etal_2016a,NeOeGaKrTh2016b}, the computational cost associated with the AD for LES/DNS computations at high Reynolds numbers and with large control horizons is still prohibitively high.

Another important issue, especially for the turbulent flow control problems, is the stability of the adjoint solvers. In general, the sensitivities that are evaluated by the adjoint solver are highly sensitive to the initial and boundary conditions at the flow regimes with high Reynolds numbers. This is actually not surprising since the turbulence phenomena, by its very nature, means that small variations in the control parameters lead to high variations in the objective function. Especially, if a wide time horizon for the flow control is chosen, the sensitivities of the specified objective function tend to an arithmetic overflow beyond a certain simulation time. The result is then an unstable adjoint solver, which evaluates sensitivity values that are not meaningful. This problem is addressed in \citep{wang13}. In contrast to URANS adjoint solvers, the usage of AD techniques does not help much alleviating this problem. The reason for this is that, the numerical scheme that is used to resolve physical instabilities that are inherent in the flow is also differentiated exactly by the AD tools, and the resulting high sensitivity values due to small scales of turbulence are simply transferred to the complete domain. Therefore, after a certain number of time iterations are performed, eventually the complete adjoint solution gets corrupted. Yet in other situations, poor numerical treatment may also lead to the same problem. Even if the physical instability is filtered out from the simulation (e.g., by using averaging or coarse grids), numerical noise introduced by an improper numerical scheme may lead to similar problems. A transient ODE example, which illustrates this problem in a detailed way can be found in \citep{Ozk16Arxiv}.   

The failure that is mentioned above can be overcome by replacing the initial value problem with the well-conditioned “least squares shadowing (LSS) problem” \citep{WANG2014, SteffiOS}. In this way, one can obtain workable sensitivities from the adjoint solver to be used in flow control problems. The drawback of the LSS method is the increased run-time and memory demand, which may be a serious problem for large-scale simulations. To decrease the run time, the LSS problem can be solved parallel-in-time \citep{Guenther2017}. As an alternative to the LSS method, the receding horizon approach \citep{Marinc2012} can be taken. In this approach the long control interval is divided into smaller sub-intervals, in which each sub-interval is small enough such that the sub-optimization problem stays controllable. In this way, an optimal control problem is treated as a group of sub-optimal problems, in which solutions to them can be achieved with classical nonlinear optimization methods. The receding horizon algorithm does not suffer from high computational cost, as the LSS method does, but this advantage comes at the expense of accuracy. 

In the present work, we aim to make a comprehensive study of the adjoint-based turbulent flow control and the associated stability issues. Thereby, we focus on the pure flow control problem of plane and round jets for noise reduction without regularization techniques that are mentioned previously. This paper is organized as follows: In Section 2, we present briefly the governing equations and the numerical method chosen for the present work. Information about the test case configurations and implementation details for the discrete adjoint solver has been also provided. The validation results for the discrete adjoint solver has been presented in Section 3. In Section 4, we present the optimization results achieved for noise reduction problem using different configurations. In Section 5, we shortly introduce the methodology used to generate the exact tangent-linear solver that has been used to study the behavior of the control sensitivities for long time horizons. In Section 6, the results of a sensitivity study obtained with the tangent-linear method are presented. Finally, some conclusions are drawn in Section 7.      

\section{Governing Equations and Numerical method}
\label{sec:numerical_method}
In the present work, the $3$D compressible Navier-Stokes equations are solved on a Cartesian 
grid to provide the primal solution
\begin{eqnarray}
\frac{\partial \rho}{\partial t} & = & \frac{\partial m_i}{\partial x_i}  \label{eqNS1}\\
\frac{\partial m_i}{\partial t} & = & -\frac{\partial p}{\partial x_i} -\frac{\partial}{\partial x_j} \rho u_j u_i + \frac{\partial}{\partial x_j} \tau_{j i}  \label{eqNS234}\\
\frac{\partial p}{\partial t} & = & -\frac{\partial}{\partial x_i} p u_i + \frac{\partial}{\partial x_i} \lambda (\gamma - 1) \frac{\partial}{\partial x_i} T\nonumber\\
 & & - (\gamma - 1 ) p \frac{\partial}{\partial x_i} u_i + (\gamma - 1 ) \tau_{i j} \frac{\partial}{\partial x_j} u_i, \label{eqNS5}
\end{eqnarray}
where $\rho$ is the density, $u_i$ is the $i$th component of the velocity vector $u$, $\gamma$ is the ratio of the specific heats, $T$ is the temperature and $\lambda$ is the heat conductivity. Furthermore, the mass flux in the $i$th direction is denoted by $m_i  =  \rho u_i$ and the viscous stress tensor $\tau_{ij}$ is given by 
\begin{eqnarray}
\tau_{i j} = \mu s_{i j} & = & \mu\left( \frac{\partial u_i}{\partial x_j} + \frac{\partial u_i}{\partial x_j} - \delta_{ij} \frac{2}{3} \frac{\partial u_k}{\partial x_k} \right),
\vspace{-.1cm}
\end{eqnarray}
where $\mu$ is the viscosity.

The above equations are discretized in space by using an optimized explicit dispersion-relation-preserving summation-by-parts (DRP-SBP) finite-difference scheme of sixth-order. As the time discretization, a two step low-dispersion-dissipation fourth-order Runge-Kutta (RK) scheme as in \citep{Hu:1996} is used. Furthermore, characteristic boundary conditions (CBC) as proposed in \citep{Lodato2008} are used to simulate the open boundaries in the jets. For the isotropic turbulence simulations periodic boundary conditions (BC) are used. Additionally, sponge regions together with grid stretching and spatial filtering are utilized to reduce the reflections near to non-periodic boundaries \citep{Foysi:2010a}. The flow field is filtered in every second RK iteration using a $10$th-order accurate low-pass filter. The filtering improves the numerical stability and also serves as a sub-grid scale model, which is equivalent to the approximate deconvolution approach to LES \citep{Stolz:1999,Mathew:2003,Mathew:2006}.

\begin{table}[h!]
\centering
{\small
\begin{tabular}{|c|cccccccccc|}
\hline
Case     & $L_x$ & $L_y$ & $L_z$ & $n_x$ & $n_y$ & $n_z$ & $\Delta_{x,min}$ & $\Delta_{y,min}$ & $\Delta_{z,min}$ & $\Delta t$ \\\hline
\jettwo  & 30    & -     & 34    & 512   & 1     & 640   & $0.04$           & -                & $0.036$          & $0.017$    \\\hline
\jetELES & 37    & 9     & 28    & 416   & 64    & 320   & $0.071$          & $0.14$           & $0.065$          & $0.03 $    \\\hline
\jetLES  & 37    & 9     & 28    & 512   & 160   & 400   & $0.051$          & $0.056$          & $0.051$          & $0.021$    \\\hline
\jetDNS  & 37    & 9     & 28    & 800   & 288   & 600   & $0.029$          & $0.031$          & $0.028$          & $0.012$    \\\hline
\end{tabular}}
\caption{Parameters of the plane jet simulations.}
\label{tab:parameterplanejet}
\end{table}

\begin{table}[h!]
\centering
\begin{tabular}{|c| c c c c c c c c c c |}
\hline
Case   & $L_x\!$ & $L_y\!$ & $L_z\!$ & $n_x$ & $n_y$ & $n_z$ & $\Delta_{x,min}$ & $\Delta_{y,min}$
 & $\Delta_{z,min}$ & $\Delta t$\\
\hline
jetG1  & 31 & 16 & 23 & 352 & 160 & 224  & 0.0675 & 0.0695 & 0.067 & 0.0286\\
jetG2  & 31 & 16 & 23 & 448 & 216 & 288  & 0.0505 & 0.049 & 0.049 & 0.0207\\
jetG3  & 31 & 16 & 23 & 640 & 288 & 384  & 0.033  & 0.033 & 0.032 & 0.0136\\
jetG4  & 31 & 16 & 23 & 1152& 512 & 640  & 0.016  & 0.015 & 0.015 & 0.0061\\
\hline
\end{tabular}
\caption{The parameters of the round jet simulations.}
\label{tab:parameterroundjet}
\end{table}

\begin{table}[h!]
\centering
\begin{tabular}{|c| c c c c c c c c c c |}
\hline
Case   & $L_x\!$ & $L_y\!$ & $L_z\!$ & $n_x$ & $n_y$ & $n_z$ & $\Delta_{x,min}$ & $\Delta_{y,min}$
 & $\Delta_{z,min}$ & $\Delta t$\\
\hline
iso64  & 1 & 1 & 1 & 64 & 64 & 64  & 0.0156 & 0.0156 & 0.0156 & 0.00219\\
iso128  & 1 & 1 & 1 & 128 & 128 & 128  & 0.0078 & 0.0078 & 0.0078 & 0.00098\\
iso256  & 1 & 1 & 1 & 256 & 256 & 256  & 0.0039  & 0.0039 & 0.0039 & 0.00013\\
\hline
\end{tabular}
\caption{The parameters of the isotropic turbulence simulations (rounded).}
\label{tab:parameteriso}
\end{table}

An important problem in industrial applications concerns the sound emission of subsonic plane and round jets and its control to suppress the radiated sound. This canonical flow problem represents a configuration, which is still simple enough allowing researchers to concentrate on the relevant physical mechanisms associated with shear flows and turbulence, without dealing with other complex effects like chemical reactions, multiple-phases, complex geometries, etc. Since the primary focus of the present work is the adjoint based optimization, the jet simulations serve as a framework for the assessment of the different adjoint approaches. In addition, we also performed forced isotropic turbulence simulations, which are good for comparing quantities like Lyapunov exponents.

Tables \ref{tab:parameterplanejet}, and \ref{tab:parameterroundjet} give an overview of the numerical parameters of the plane and round jet simulations that are performed in this work.
For the plane jet simulations, the domain lengths 
$L_i$ are normalized by the jet diameter $D$. The Reynolds number, based on the diameter is set to $Re =U_j \rho_j D/\mu_j=2000$ and the Mach number is set to  $Ma=U_j/c_j=0.9$ using centerline values for all simulations. The number of grid points are represented by $n_i$ for the different directions, respectively, with associated minimum grid-spacing $\Delta_{i,min}$. 
$\Delta t$ indicates the time-step used during the optimization computations non-dimensionalized by $D/U_j$. The subscript $j$ denotes mean values at the jet inflow plane. Similarly, for the round jet simulations $L_i$ is the length of the computational domain in $i$th. direction. For the jets the reference length $D_j$ and reference velocity $U_j$ are chosen to be the jet diameter and velocity at the inflow. In Table \ref{tab:parameteriso}, numerical parameters used for the forced isotropic turbulence simulations are given. For the isotropic turbulence, the extent of the computational domain $L_i$ and the rms-velocity over the whole computational domain serve as reference values for non-dimensionalization. The Mach number $Ma_{rms}$ is set as $0.2$ for the isotropic turbulence. In order to sustain the turbulence the flow has to be forced. In the jet simulations the flow is forced at the inflow by explicitly setting the inflow part of the CBC. Precurser simulations, which are fed into the domain to provide realistically correlated inflow data can be found in \citep{Foysi:2010a,Marinc2012}.
For the case of isotropic turbulence the forcing method proposed in \citep{Petersen} was implemented. In this work, only the solenoidal part of the thirteen modes with the lowest wave length were perturbed. Note that this forcing term allows us to choose explicitly the dissipation rate $\varepsilon$, therefore enabling an easy evaluation of the Kolmogorov length scale $\eta = (\nu^3/\varepsilon)^{1/4}$.

The maximal Lyapunov exponent (MLE) can be computed by first solving the linearized Navier-Stokes equations with some arbitrary initial perturbation. As soon as the solution exhibits exponential growth in $L2$ norm of the flow quantities, the MLE is obtained by a fit to this exponential following the procedure described in \cite{Kuptsov}.\\

\begin{figure}
  \begin{center}
    \includegraphics[scale=0.45,trim={0.5cm 0.3cm 0cm 0.3cm},clip]{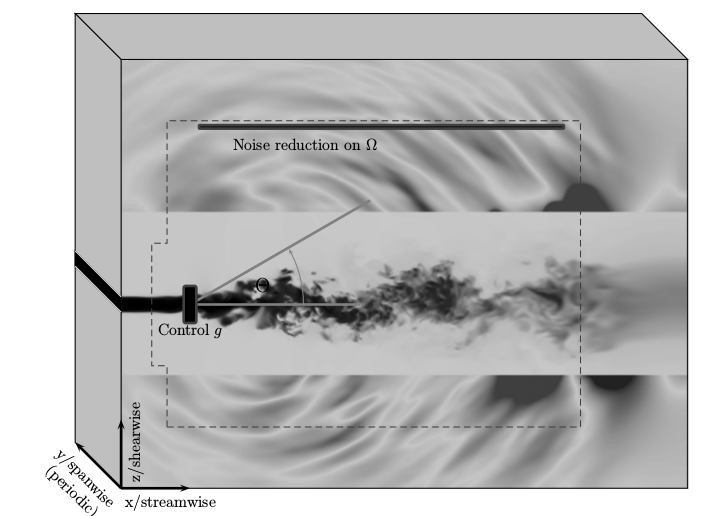} 
    \caption{An illustration of the control setup.}
    \label{fig:control_setup}
  \end{center}
  \end{figure}

As already mentioned, the investigation of the adjoint approaches is centered on the problem of reducing sound emission of compressible turbulence jets. A space and time dependent heating/cooling source term in an area near the jet inflow acts as the control. The noise emitted by the jet is measured by a cost functional, which is defined as an integral over the square of the pressure fluctuations in the far field
\begin{equation}
\Im = \int_{\Omega} \int_T \left(p(\mathbf{x},t)-\overline{p}(\mathbf{x}) \right)^2 dt~d\Omega,
\label{eqCost}
\end{equation}
where $T$ is the control interval length in time, $\Omega$ is a small volume in the far-field of the jet and $\bar p$ denotes the temporal average of the pressure over the interval $T$. The control setup is illustrated in Fig. \ref{fig:control_setup}. Optimization studies using configurations similar to this can be found in \cite{Kim:2010,Wei:2006}.

As far as the gradient evaluation method is considered, the adjoint method has been chosen for its computational efficiency. Based on the flow solver, which has been introduced previously, both discrete and continuous adjoint versions have been developed. The implementation details of the continuous adjoint can be found in \cite{Marinc2012} and is not repeated in this paper for the sake of brevity. On the other hand, the hand-discrete adjoint solver includes some subtle details. Therefore, the key implementation issues regarding the development of the discrete adjoint solver is briefly introduced in the following.

Within a time iteration, $N$ Runge-Kutta sub-steps together with filtering, control and initial conditions can be written generally as
\begin{eqnarray}
\bold k_0 &=& \bold 0 \\
\bold u_0 &=& \bold u_{init}\\
\bold k_s &=& \alpha_{s-1}\bold k_{s-1} + \Delta t \left[
\bold R(\bold u_{s-1}) + \sum_{i=0}^M \gamma_{s-1,i}\mathbf\Phi_i\right]\\
\bold u_s &=& F_s[\bold u_{s-1} + \beta_{s-1}\bold k_s],\qquad s \in \{1,\ldots,N\},
\end{eqnarray}
where $\bold R$ is the right hand side (RHS) of the discrete Navier-Stokes equations, 
$F_s$ is the discrete representation of a filter operator at step $s$ in case of LES, 
$\mathbf\Phi_i$ are $M+1$ control vectors and $\gamma_{s,i}$ are scalars chosen such that a linear 
interpolation between the control vectors is achieved.

Given some cost functional $\Im = \sum_{s=0}^N \Im_s(\bold u_s)$
 the Lagrangian is defined as
\begin{eqnarray}
L &=&\sum^N_{s=0}\Im_{s}(\bold u_{s}) - \sum^N_{s=1}\boldsymbol{\xi}^T_s\left[
\bold k_s - \alpha_{s-1}\bold k_{s-1} - \Delta t \left\{
\bold R(\bold u_{s-1}) + \sum^M_{i=0}\gamma_{s-1,i}\mathbf\Phi_i\right\}\right]\nonumber\\
&&-\sum^N_{s=1}\boldsymbol{\omega}^T_s(\bold u_{s} - F_s[\bold u_{s-1} + \beta_{s-1}\bold k_s])
- \boldsymbol{\xi}^T_0 \bold k_0 - \boldsymbol{\omega}^T_0 [\bold u_{0} - \bold u_{init}],
\end{eqnarray}
where the vectors $\boldsymbol{\xi}_s$ and $\boldsymbol{\omega}_s$ are the Lagrangian multipliers. 
Using the variational form of this equation, transposing and rearranging the terms with respect to the variations leads to the adjoint Runge-Kutta integration scheme
\begin{eqnarray}
\boldsymbol{\omega}_N &=& \left(\left.\PD{\Im_{N}}{\bold u}\right|_{\bold u_N}\right)^T\\
\boldsymbol{\xi}_N&=& \beta_{N-1} F^T_N\boldsymbol{\omega}_N\\
\boldsymbol{\omega}_s&=&F^T_{s+1}\boldsymbol{\omega}_{s+1} + \Delta t
 \left(\left.\PD{\bold R}{\bold u}\right|_{\bold u_s}\right)^T\boldsymbol{\xi}_{s+1}+
\left(\left.\PD{\Im_{s}}{\bold u}\right|_{\bold u_s}\right)^T \; s \in \{0 ,\ldots,N-1\}\label{eq:om}\\
\boldsymbol{\xi}_s &=&\alpha_{s}\boldsymbol{\xi}_{s+1} + \beta_{s-1}F^T_s \boldsymbol{\omega}_s \quad s \in \{1,\ldots,N-1\}\\
\boldsymbol{\xi}_0 &=&\alpha_{0}\boldsymbol{\xi}_1.
\end{eqnarray}
Using  $\boldsymbol{\xi}_s$, the gradient of the cost functional is given by
 \begin{equation}
\left(\frac{d\Im}{d\mathbf\Phi_i}\right)^T =
 \left(\frac{d L}{d\mathbf\Phi_i}\right)^T = \Delta t\sum^N_{s=1}\gamma_{s-1,i}\boldsymbol{\xi}_s
,\qquad i\in \{0,\ldots,M\}.
\end{equation}
Note, that in equation (\ref{eq:om}) the transpose of the linearized Navier-Stokes operator $ \left(\left.\PD{\bold R}{\bold u}\right|_{\bold u_s}\right)^T$ appears. 
When implementing the discrete adjoint equations the hardest part is usually the implementation 
of this operator, as the components of the linearized Navier-Stokes operator are for most cases  not known explicitly and thus its transpose cannot be computed directly.
In order to compute the transpose of some linear operator $D$ (e.g. a matrix representing a 
discretized derivative operator) we utilize its expression as
\begin{equation}D =\sum^N_{i=1}D_i,\label{eq:split}\end{equation}
where $D_i$ has the entries of $D$ in its $i$th row and is zero everywhere else. 
Consequently, the operator $D_i^T$ is zero except for the $i$-th column. 
The product $D_i\bold a$, with some arbitrary vector $\bold a$, is computed by multiplying grid 
points adjacent to the $i$th grid point with the corresponding finite difference coefficients 
and evaluating the sum of these values at the $i$th grid point. In a similar fashion, 
the product $D^T_i\bold a$ can be computed by setting the corresponding adjacent grid points to the
product of the $i$th grid point with a finite difference coefficient. 
The implementation of the transposed linear operator-product is similar to the one described in \cite{Pando}. With the operator splitting in Eq. (\ref{eq:split}), the transpose of the full Navier-Stokes operator can be computed for each grid point separately as $D^T = \sum^N_{i=1} D_i^T$.
Consequently, different parts of the computational domain, e.g. boundary and inner schemes,
can be treated on its own. Furthermore, only the entries of the $i$th row of $D$ are needed 
for computing the contribution of the $i$th grid point to the transposed right hand side. Therefore,
the derivative coefficients of adjacent grid points do not need to be known, which eases the parallel
implementation. The procedure described above may be interpreted such that the matrix vector 
product is computed row-wise, while the product with the transposed matrix is computed column-wise.
The Navier-Stokes equations consist of several sums and products of linear operators, therefore
the identity
\begin{equation}
(AB + CD)^T = B^T A^T + D^T C^T
\end{equation}

for arbitrary matrices $A, B, C$ and $D$ has to be used to make the above implementation 
applicable to the cases considered in this work. Splitting all the derivative operators within the 
Navier-Stokes equations according to the preceding equation and the Eq. (\ref{eq:split}) can become a tedious task. Therefore, a program is implemented that automatically generates the source code for computing the desired matrix vector products for any given PDE.

\section{Validation of the Discrete Adjoint}
The implementation of the hand-discrete discrete adjoint can be error-prone, therefore several tests have been performed to check the validity of the approach used in this work. To ensure that the linearization of the Navier-Stokes equations is correctly implemented, 
results from the sensitivity operator were compared with the sensitivities obtained through complex 
differentiation of real functions:
\begin{align}
   s(x)&=Re(s(x+i h))+\mathcal{O}(h^2), \\
   s^{(1)}(x)&=\frac{Im(s(x+i h))}{h}+\mathcal{O}(h^3),
\end{align}
where $s$ is an arbitrary function, $i$ is the imaginary unit and $Re(\ldots)$ and $Im(\ldots)$ are the real and imaginary parts. With this formulation the sensitivities of an arbitrary function can be computed by just exchanging real values with the complex ones in the function. Furthermore, the sensitivity is calculated without performing subtractions, which means that the step size $h$ can be chosen very small without having cancellation errors.
Tests showed that the difference of the sensitivities obtained by using the direct 
implementation and with that using complex differentiation was of the order of machine precision.
Note that the RK-iteration is already linear so that making additional checks for its correct is linearization unnecessary.

The comparison between the sensitivity and discrete adjoint equations was done using the identity 

\begin{equation}
\mathbf{a}^T\left(\left.\frac{\partial \bold R}{\partial\mathbf{\Phi}}\right|_{\mathbf{\Phi_{s}}}\right)^T\mathbf{b}=\mathbf{b}^T\left(\left.\frac{\partial \bold R}{\partial\mathbf{\Phi}}\right|_{\mathbf{\Phi_{s}}}\right)\mathbf{a}
\end{equation}
for randomly chosen vectors $\mathbf{a}$ and $\mathbf{b}$, with ${\bold R}$ being the Navier-Stokes 
operator. This test revealed that, the transpose of the linearized Navier-Stokes operator is implemented correctly and is accurate up to machine precision.
The complete implementation of the full adjoint system, given some arbitrary perturbation $\mathbf{\Phi}'$ of the control, the linear response of the cost functional with respect to that perturbation can be obtained either from the solution of the sensitivity equations (the linearized state equations), or by the first order term of a Taylor series expansion via the discrete adjoint:
\begin{equation}
\underbrace{\frac{\partial L}{\partial\mathbf{\Phi}}}_{\text{gradient}}\mathbf{\Phi}'
=\frac{d\Im}{d\mathbf{\Phi}}\mathbf{\Phi}'=
\frac{d\Im}{d\mathbf{u}}\frac{d\mathbf{u}}{d\mathbf{\Phi}}\mathbf{\Phi}'+\PD{\Im}{\mathbf{\Phi}}\mathbf{\Phi'}=
   \frac{d\Im}{d\mathbf{u}}\underbrace{\mathbf{u}'}_{\text{sensitivity}}+\PD{\Im}{\mathbf{\Phi}}\mathbf{\Phi'}
   \label{eqSensGradIdentity}
\end{equation}
where $L$ is the Lagrangian of the system at hand and the superscript $(\cdot)'$ denotes the variation with respect to the control. 
The identity in eq. \eqref{eqSensGradIdentity} involves integration over the full control horizon and was used to validate the implementation of the discrete adjoint further. 
Several tests with small Gaussian shaped perturbations with a finite temporal support were performed for the cases listed in \cite{Marinc2012}. 
A selection of these tests are depicted in Table \ref{tab:validate} and it becomes apparent that the both sides of the eq. \eqref{eqSensGradIdentity} match perfectly. 
Although the RHS is transposed accurately up to machine precision, the full discrete adjoint over 
the complete time horizon is slightly less accurate. A possible reason are the round-up
errors that occur during the RK iterations as well as the integration of long time horizons. 
Furthermore, the discrete adjoint is utilized as a reference solution  for validation of various
other adjoint formulations. 

\begin{table}[h!]
\centering
\begin{tabular}{|c|c c c c|}
\hline
case& pos.& \#RK& $\PD{L}{\mathbf{\Phi}}\mathbf{\Phi}'$ & $\PD{\Im}{\boldsymbol{u}} 
{\boldsymbol{u}'} + \PD{\Im}{\mathbf{\Phi}}{\mathbf{\Phi}}'$\\
\hline
DNS2D &100& 10000   & -2.19978130134220$\cdot10^{-4}$ & -2.19978130134229$\cdot10^{-4}$\\
      & 900 & 10000 & -1.44380102098536$\cdot10^{-5}$ & -1.44380102098535$\cdot10^{-5}$\\
      &1700 & 10000 & -1.17132271501490$\cdot10^{-5}$ & -1.17132271501492$\cdot10^{-5}$\\
\hline
ELES3D &100 & 2400 & -5.37171660047461$\cdot10^{-6}$ & -5.37171660047416$\cdot10^{-6}$\\
       &200 & 2400 & -6.73097452538381$\cdot10^{-6}$ & -6.73097452538368$\cdot10^{-6}$\\
\hline
LES3D &100 & 6400 & -2.25323809359032$\cdot10^{-4}$ & -2.25323809359047$\cdot10^{-4}$\\
      &200 & 6400 & -1.70187568261517$\cdot10^{-4}$ & -1.70187568261540$\cdot10^{-4}$\\
\hline
DNS3D &100 & 7000 & -3.11011806150672$\cdot10^{-4}$ & -3.11011806150674$\cdot10^{-4}$\\
      &400 & 7000 & -2.85005052892766$\cdot10^{-5}$ & -2.85005052892798$\cdot10^{-5}$\\
\hline
\end{tabular}
\caption{Comparisons of both sides of eq. (\ref{eqSensGradIdentity}) for different perturbations 
and for different jet simulations. The perturbation was initiated after $pos.$ Runge-Kutta
iterations from the beginning of the control interval, extending over $\#RK$ Runge-Kutta iterations. 
The maximum normalized amplitude of the perturbations was $10^{−4}$.}
\label{tab:validate}
\end{table}
The gradient is used in a low storage Limited memory Broyden-Fletcher–Goldfarb-Shanno (L-BFGS) optimization scheme together with a line search method using the Wolfe condition.

\section{Optimization Results with Discrete and Continuous Adjoint Solvers}

\begin{figure}
\begin{center}
\includegraphics[scale=0.4]{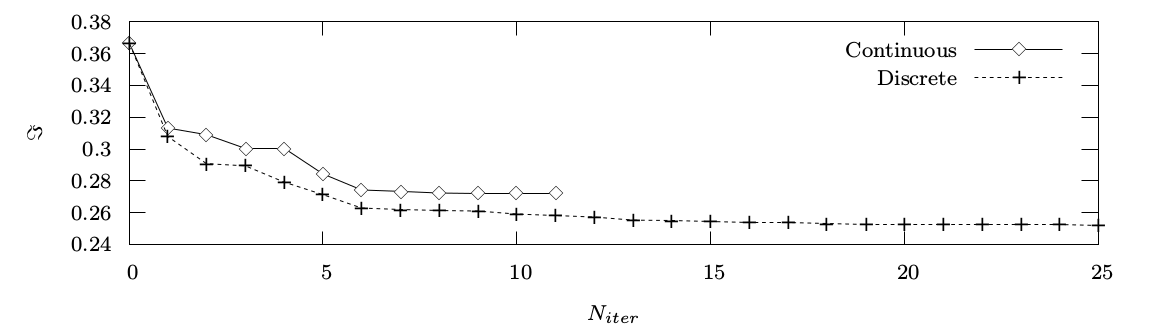} 
\caption{Comparison of the cost functional over L-BFGS-iterations for case \jettwo\ involving a long control interval. The discrete optimization clearly outperforms the continuous optimization.}
\label{figJet2DoptLongIter}
\end{center}
\end{figure}

\begin{figure}
\begin{center}
\includegraphics[scale=0.4]{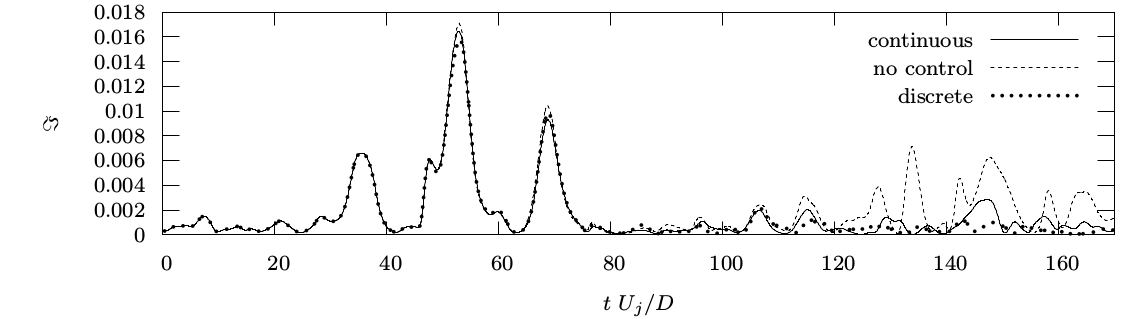} 
\caption{Comparison of the cost functional over time for the \jettwo\ case with a long control interval.}
\label{figJet2DoptCostTime}
\end{center}
\end{figure}

Figure \ref{figJet2DoptLongIter} shows the reduction of the cost functional for the \jettwo\ case for both continuous and discrete optimization over the number of L-BFGS iterations. The same behavior over time can be observed in Figure \ref{figJet2DoptCostTime}, which also clearly indicates the time horizon affected by the control.
It is clearly seen, that the discrete adjoint optimization is able to reduce the sound pressure level (SPL) more than the continuous approach. However, in both cases it was not possible to effectively influence the SPL below $t U_j/D=100$.
It should be noted that this interval could be controlled successfully when optimizing over this shorter intervals. This was already observed in \cite{Marinc2012} while dealing with the continuous approach. 
In this regard, the discrete approach was considered to be more promising, since the differences in the numerical treatment and the boundary conditions used while discretizing the continuous adjoint equations may lead to inconsistencies compared to the primal flow equations. These differences result that the discretized adjoint equations are not the exact dual counterpart to the discretized primal equations.
The fact that long control intervals are difficult to control even with accurate gradient information suggests that dividing the optimization problem into a set of sub-optimal problems with shorter time horizons might be more efficient. Such strategies are realized using the receding-horizon approach, as done and proposed in \cite{Collis:2000,Marinc2012}. 

\begin{figure}
\begin{center}
\includegraphics[scale=0.4]{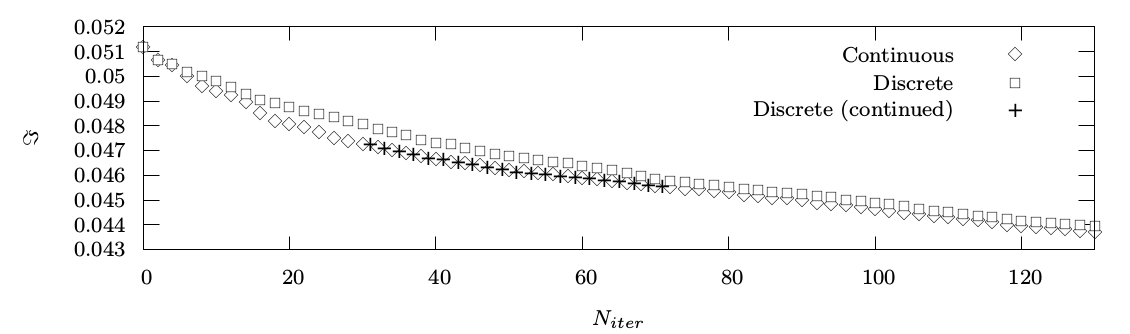} 
\caption{Cost functional over L-BFGS iterations for the ELES3D case with a short control interval.}
\label{figJetELESoptIter}
\end{center}
\end{figure}

When applying the procedures to the $3$D simulations, a stronger reduction in the sound-pressure levels of the $2$D simulations has been observed. Another observation is that, the continuous optimization performed minimally better in terms of cost functional reduction in this case. To check whether this behavior was triggered by the choice of initial condition, the optimization using the discrete adjoint has been performed using the restart solution obtained from continuous adjoint optimization with $30$ iterations. 
The result of this optimization study as well as pure discrete and continuous adjoint optimizations are shown in Figure \ref{figJetELESoptIter}. In both cases, the cost functional can be significantly reduced, and only a negligible difference between both approaches has been observed. This behavior supports the conjecture, that the differences observed in Figure \ref{figJetELESoptIter} are due to the high dimensionality of the control vector, and thus the 
increased complexity of the cost functional response surface.
For short control horizons, as used in this case, the performance of the optimization is similar for the continuous and discrete approaches as both the approaches still calculate accurate gradients. In both cases, approximately $3$dB reduction in the overall SPL has been achieved.

\begin{figure}
\begin{center}
\includegraphics[scale=0.4]{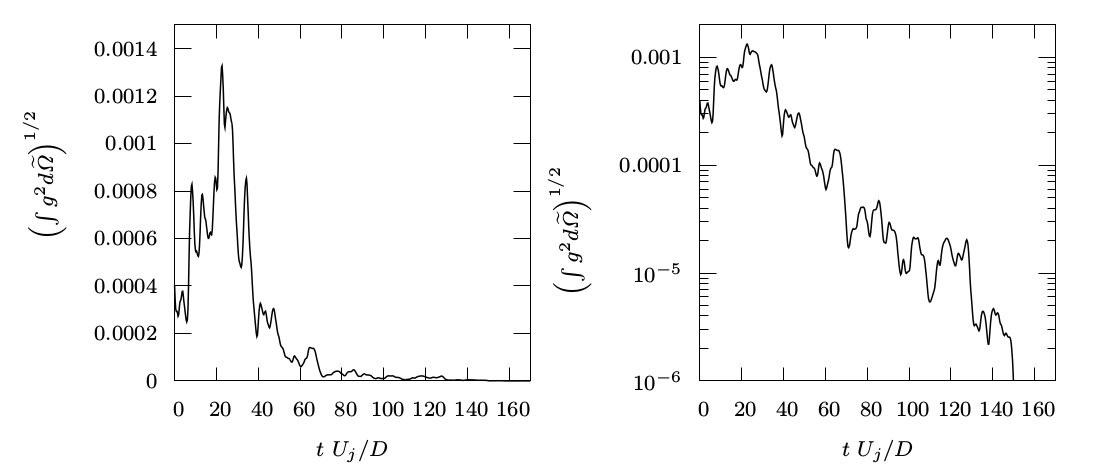} 
\caption{Left: The quantity $\int g^2 d \tilde{\Omega}$ as a measure for the gradient-”strength” over time, where $g$ is the gradient of the cost functional and $\tilde{\Omega}$ is the whole computational domain. Right: The same quantity in log-scale.}
\label{figJet2DabsgradVLong}
\end{center}
\end{figure}

To illustrate the problems emerging from the long control horizons in more detail, Figure \ref{figJet2DabsgradVLong} shows the ''energy`` of the gradient obtained with a control vector set to zero.
It can be observed that the norm of the gradient decreases strongly with time, which is a hint that the problem is very ill-conditioned and is thus a source for optimization inefficiency.
Furthermore, as the control is basically a linear superposition of gradients obtained during the optimization iterations, the control has high amplitudes only at the beginning of the control interval but it is negligible for the rest. Consequently, a significant part of the control interval is practically unused.

\begin{figure}
\begin{center}
\includegraphics[scale=0.4]{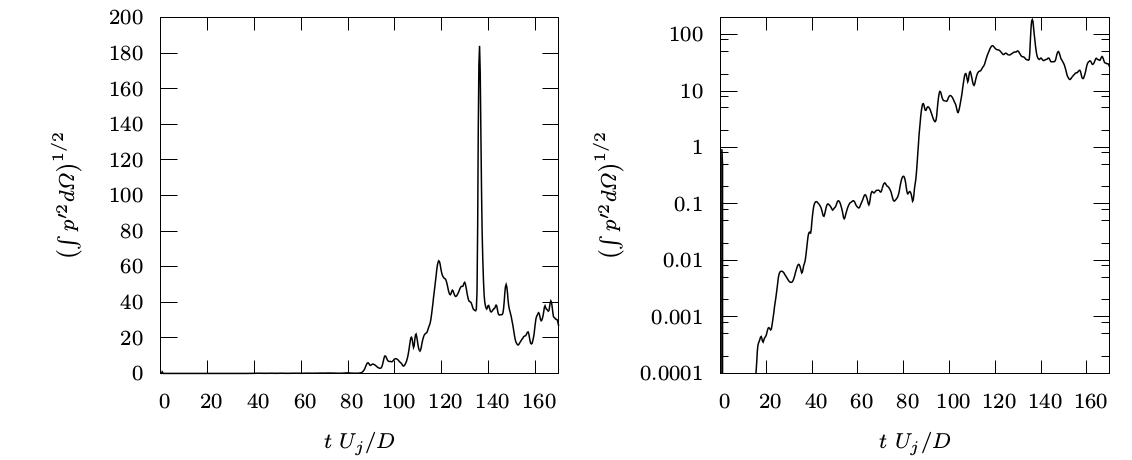} 
\caption{Left: The quantity $\int p^{\prime2} d\Omega$ as a measure for the linear response of the cost functional to a control perturbation. It can be seen that only times
$t \geq 100$ are influenced noticeably by the control. Right: The same quantity in log-scale.}
\label{figJet2DCostSensVLong}
\end{center}
\end{figure}

Another problem with the long control horizons becomes apparent by investigating Figure \ref{figJet2DCostSensVLong}, which shows the quantity $\int p'^2 d\Omega$ as a measure of the strength of the linear response of the cost functional due to a control perturbation. The gradient computed in Figure \ref{figJet2DabsgradVLong} was chosen as the control perturbation ($\mathbf g:= \mathbf\Phi'=\left(\frac{d\Im}{d\mathbf{\Phi}}\right)^T$) and $p'$ constitutes the linear response of the pressure to this perturbation.
It can be observed that the linear response of the cost functional increases exponentially with time and large fluctuations can occur, leading to very large values of the gradient sensitivities in some cases.
As the gradient used for optimization contains no reliable information for the non-linear regime, the optimization scheme will in most cases select a step size, in which the linear terms cannot dominate. Thus, the linear response of the cost functional to a perturbation gives an estimate of the change of the cost functional over one optimization iteration. This implies, that due to the strong increase of the amplitude of the sensitivities, they become often hardly differentiable, only a short interval at the end of the simulation is actually controlled.
It should be noted that the strong amplitude growth is not related to numerical instabilities but rather to instabilities that are inherent to the linearized Navier-Stokes equations, acting on the initial and boundary conditions.

\begin{figure}
\begin{center}
\includegraphics[scale=0.4]{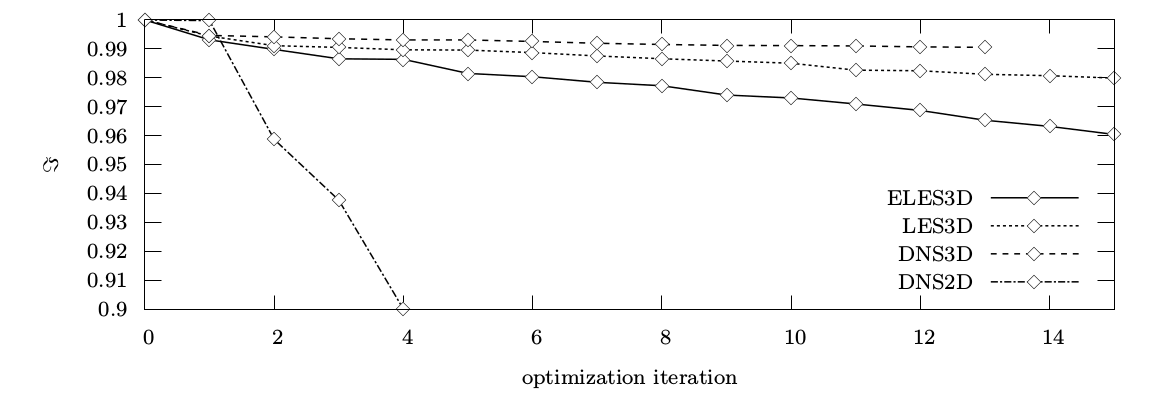} 
\caption{Different optimization runs for the cases \jettwo, \jetELES, \jetLES\ and \jetDNS\ }
\label{figJetsInCompRes}
\end{center}
\end{figure}

Reducing the dimensionality of the non-linear system improves the situation, thus possibly giving the impression of successful optimization over given control horizons for selected time-dependent
problems. This is illustrated using $2$D and $3$D simulations, including lower resolution LES, by artificially decreasing the control space dimension. Cases \jettwo, \jetELES, \jetLES\ and \jetDNS\ were optimized using the L-BFGS scheme with a control interval length of $T=70$. The cost function reduction relative to the value without control is compared for these cases in Figure \ref{figJetsInCompRes}. For a better comparison the cost function values are normalized with their initial values in each test.
It becomes obvious that the efficiency of the optimization decreases with the increasing resolution of the $3$D simulations.
This is reasonable as the range of scales resolved in the flow solution increases with the resolution, in addition to the increase in the control space dimension. This is due to the fact that the volume of the controlled area was kept fixed in order to make the different cases physically comparable, thereby increasing the control space dimension with increasing resolution. These smaller scales likely have an adverse effect on the controllability of the system.

The Courant-Friedrichs-Lewy (CFL) criterion ensures additionally, that the time-step decreases with the increasing resolution, and consequently 
the number of RK iterations is required to be higher. Therefore, the dimension of the control vector also increases, since the cost functional involves an integration over the time horizon. In a second experiment, the control space dimension was reduced by using interpolation between selected time steps. 
An interpolation point for the control was set every second RK iteration for the cases \jettwo, \jetELES\ and \jetLES\ and every third RK iteration for the case \jetDNS.
Thus, the number of interpolation points in time increases with the resolution, also leading to an increased control space dimension.
The number of control variables can be reduced by choosing an independent control parameter at only every $n^{th}$ grid point.
The missing control function values in between are obtained using interpolation.
In the following, this kind of control is referred to as \gapCtrl{$\alpha\beta\gamma$}{$\delta$}, where $\alpha$, $\beta$ and $\gamma$ gives the number of grid points between the sampling points for the control in three spatial directions and $\delta$ means that control values are given in every $\delta$th RK iteration.
For example, \gapCtrl{111}{1} means that the control is active at every grid point in the controlled region and at every RK iteration.
The interpolation in the spatial directions is achieved using a Catmull-Rom spline \citep{Marschner:1994}.
Controlling only every $n^{th}$ grid point directly means roughly that only wave lengths below an $n^{th}$ of the Nyquist wave lengths can be controlled directly.
Consequently, by choosing the gap size between control supporting points a trade off is made between controllability and the reduction in control space dimension.
It should be noted, that the highest wave numbers are not captured by the FD derivative operators 
or cut off using filtering.
Therefore, wave lengths near the Nyquist frequency can not be efficiently controlled, which reduces the frequency band effectively controllable by scheme $gap111$ below the Nyquist frequency.

\begin{figure}
\begin{center}
\includegraphics[scale=0.4]{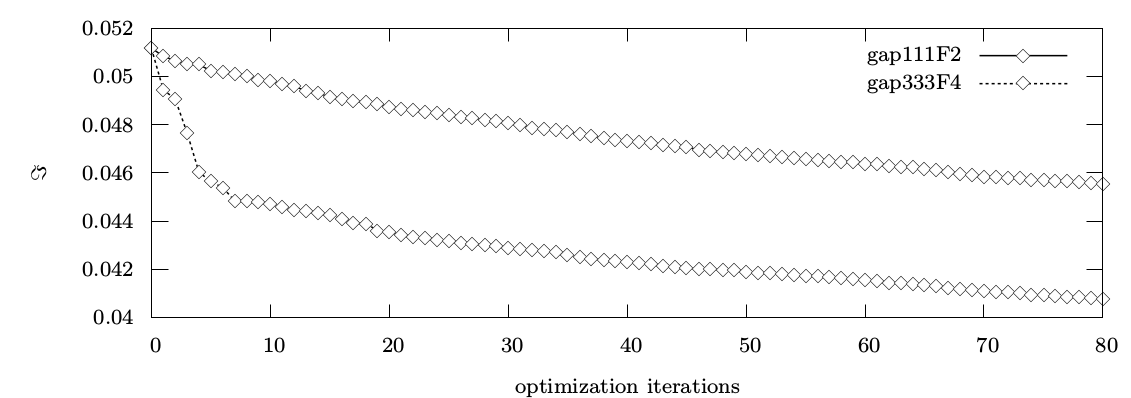} 
\caption{Optimization histories for the \jetELES\ case}
\label{figJetsInCompGapELES3D}
\end{center}
\end{figure}

\begin{figure}
\begin{center}
\includegraphics[scale=0.4]{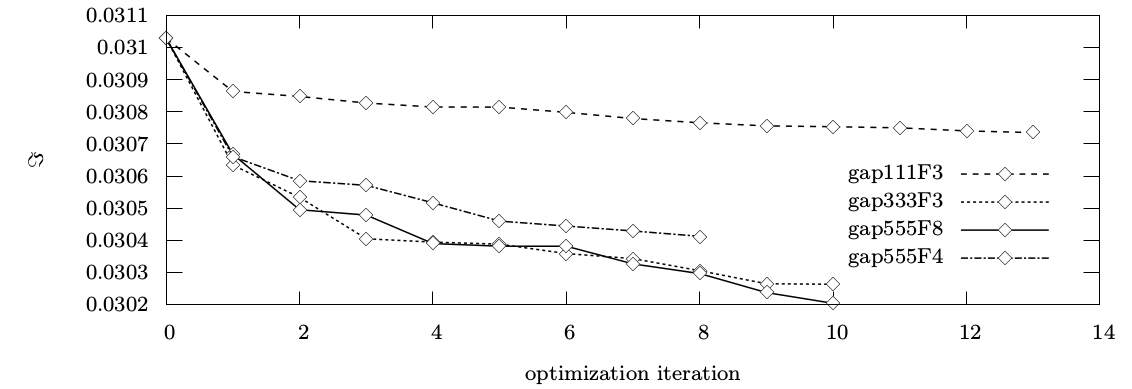}  
\caption{Optimization histories the \jetDNS\ case}
\label{figJetsInCompGapDNS3D}
\end{center}
\end{figure}

Figures \ref{figJetsInCompGapELES3D} and \ref{figJetsInCompGapDNS3D} show the history of the cost functional over optimization iterations for different gaps for the cases \jetELES\ and \jetDNS\ using a control interval length of $T=70~D/U_j$. One should have in mind that for case \jetDNS\ due to the computational cost only a moderate number of optimization iterations with only one initial condition were performed.
Thus, the specific behavior observed in Figure \ref{figJetsInCompGapDNS3D} should be interpreted with care.
However, it can be clearly observed that the efficiency of the optimization could be increased successfully by decreasing the control space dimension for both cases.
It can be observed in Figure \ref{figJetsInCompGapDNS3D} that cases \gapCtrl{333}{3}\ and \gapCtrl{555}{8}\ perform comparably well, whereas the reduction is lower for case \gapCtrl{555}{4}.
This indicates that there is not a monotonic relation between optimization efficiency increase and control space reduction.

For the case of a DNS all physically relevant effects are resolved and one would expect that the linear instabilities of the system are solely determined by the physical properties of the system.
For an LES, however, smaller scales are excluded from the numerical simulation. These scales can effect the maximum Lyapunov exponent, which we thus expect to depend on the numerical parameters chosen for LES.

\begin{figure}
\begin{center}
\includegraphics[scale=0.4]{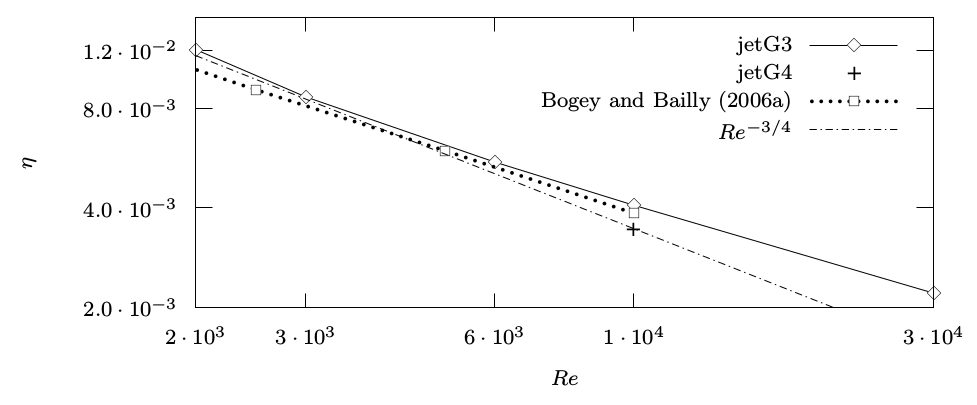}   
\caption{Kolmogorov length estimations for round jets with different Reynolds numbers.}
\label{figRkolmogorovLRe}
\end{center}
\end{figure}

To quantify the smallest physically relevant scales, Figure \ref{figRkolmogorovLRe} shows estimates of the Kolmogorov length for different Reynolds numbers, obtained from computations performed with grids \gc\ and \gd.
The estimates have been calculated via the relations given in \cite{Pope}:

\begin{align}
   \eta &= \left(\frac{\nu^3}{\epsilon}\right)^{1/4} \label{eqEstimationKolmogorov}\\
   \epsilon &=\nu\overline{\frac{\partial \widehat{u}_i}{\partial x_j}\frac{\partial \widehat{u}_i}{\partial x_j}+\frac{\partial \widehat{u}_i}{\partial x_j}\frac{\partial \widehat{u}_j}{\partial x_i}} \\
   \nu &= \frac{\mu(T_{ref})}{\rho_{ref}},
\end{align}
where $\overline{\cdots}$ denote the Reynolds average and $\widehat{u}=u-\overline{u}$.
It should be noted that Eq. \eqref{eqEstimationKolmogorov} is derived from the incompressible Navier-Stokes equations.
However, as only unheated jets are considered in this work, the density fluctuations are expected to be small enough for Eq. \eqref{eqEstimationKolmogorov} to give a useful approximation.
One can find a reasonable agreement with the lengths estimations given in \cite{Bogey:2006}.

\begin{figure}
\begin{center}
\includegraphics[scale=0.4]{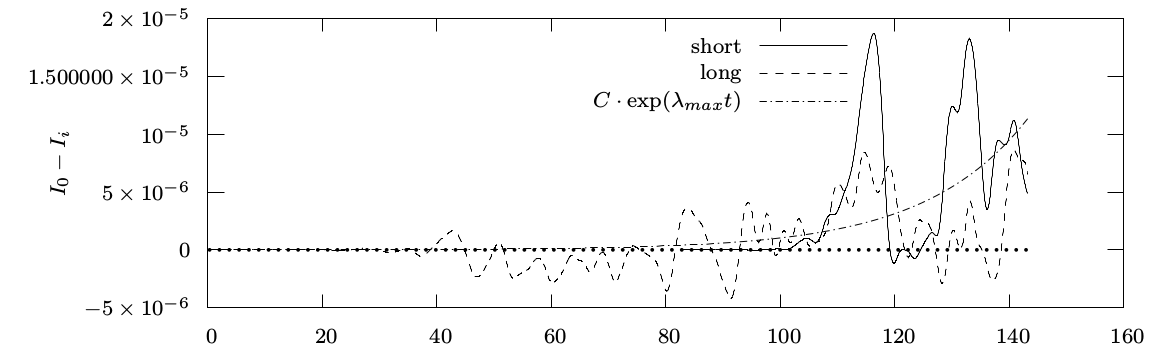}    
\caption{Reduction of the cost function over time for control intervals with length $T = 60$ and $T = 143$ and the grid jetG1.}
\label{figR1optCalcCost}
\end{center}
\end{figure}

To illustrate the impact of the MLE on the optimization, the case \ga\ has been optimized on two intervals with different temporal lengths of $T=60$ and $T=143$.
Figure \ref{figR1optCalcCost} shows the reduction of the cost function over time after eight optimization iterations for these two interval lengths.
Again it is clear that only a final interval of length $\Delta T\approx 40$ is successfully 
controlled. 
The figure also shows an exponential function growing with the rate of the MLE.
It can be seen that the length of the effectively controlled interval corresponds roughly with the average linear response of the cost functional estimated by this MLE.

\begin{figure}
\begin{center}
\includegraphics[scale=0.4]{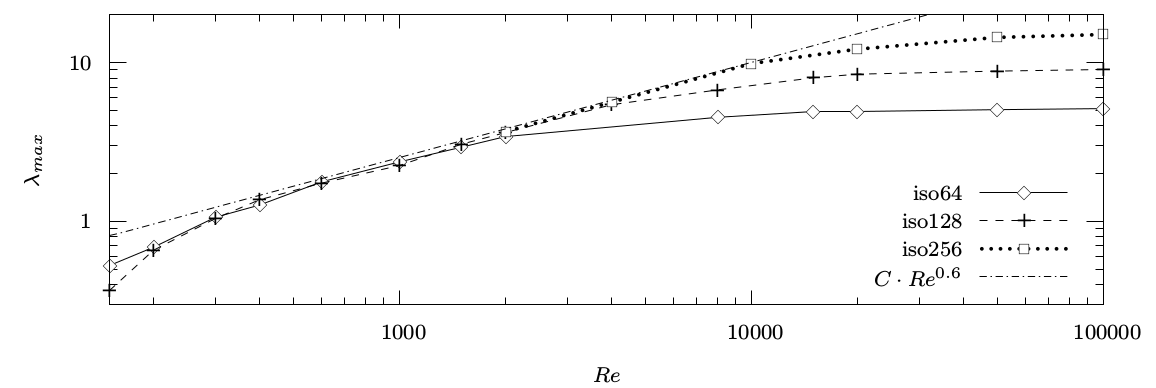}   
\caption{Lyapunov exponent with maximal real part versus Reynolds number for forced isotropic turbulence with different resolutions and Reynolds numbers.}
\label{figIsoLyapunovRe}
\end{center}
\end{figure}

In the following the maximum Lyapunov exponent is determined as a function of Reynolds number.
Due to the easier access to theory, first forced isotropic turbulence case is considered, in which
the dissipation rate $\epsilon$ is fixed.
It has been reasoned that the MLE is proportional to the smallest time scale present in the computation \citep{Crisanti:1993}. Using the relations $\tau=\sqrt{\frac{\nu}{\epsilon}}$ and $\nu=\frac{1}{Re}$ (non-dimensionalized) one 
obtains the estimate $\tau\propto \sqrt{\frac{1}{Re}}$, where $\tau$ is the Kolmogorov time scale, the smallest physical time scale in the flows considered.
In Figure \ref{figIsoLyapunovRe}, the expected behavior can be roughly seen for smaller Reynolds numbers ($Re\lesssim 4000$), where one observes a scaling with $Re^{0.6}$. Note that the above derivation does not take compressibility effects into account, which may result in the observed difference in the scaling behavior.
For larger Reynolds numbers the MLE reaches a plateau as soon as resolution is insufficient to
resolve the required length and time scales, limiting the MLE.

\begin{figure}
\begin{center}
\includegraphics[scale=0.4]{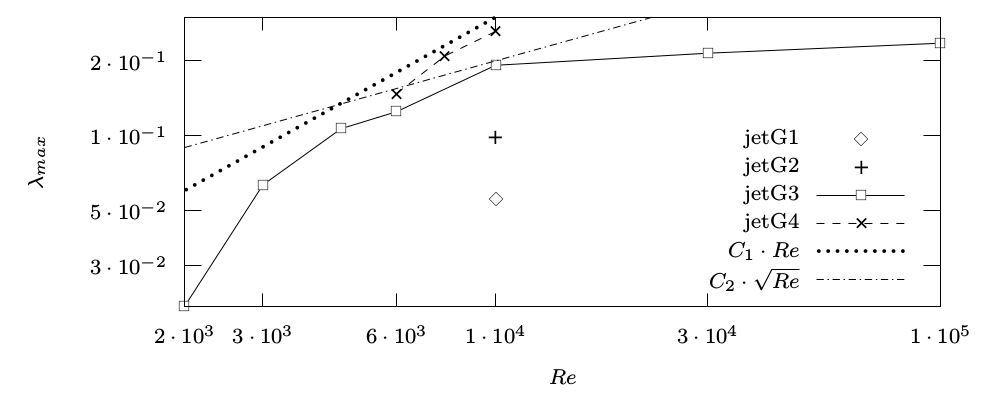} 
\caption{Lyapunov exponent with maximal real part versus Reynolds number for the round jet with different
resolutions.}
\label{figR3LyapunovRe}
\end{center}
\end{figure}

Figure \ref{figR3LyapunovRe} shows the MLE, now determined for the round jet calculations
providing information for different resolutions and Reynolds numbers. 
The figure shows the effect of increasing the resolution for a fixed Reynolds number of $Re=10000$.
As is to be expected the MLE increases with increasing resolution, supporting the behavior
seen in the optimization calculations above (Figure \ref{figJetsInCompRes}).
As already found for the forced isotropic turbulence the increase of the MLE with increasing Reynolds number reaches a saturation for very large Reynolds numbers, where the smallest scales are not sufficiently resolved. Still, both cases jetG3 and jetG4 considered for this investigation indicate an almost linear scaling with Reynolds number itself.
Opposite to the isotropic turbulence case, where intrinsic compressibility effects are dominant,
these jet cases have strongly varying fluid properties, changing the Kolmogorov length and time scales. 
As a consequence the increase of the MLE is larger than estimated above for the isotropic turbulence
simulations.

\begin{table}
	\centering
		\begin{tabular}{|c|c|}
			\hline
			LES-model & max. Lyapunov exponent \\\hline
			expl. filter ($10^{th}$ order) & $9.85\cdot 10^{-2}$ \\
			DMM1 from \cite{Martin:2000}   & $3.75\cdot 10^{-2}$ \\
			Smagorinsky (C=0.01)           & $3.12\cdot 10^{-2}$ \\
			dynamic Smagorinsky            & $1.78\cdot 10^{-2}$ \\
			Smagorinsky (C=0.02)           & $1.45\cdot 10^{-2}$ \\\hline
		\end{tabular}
	\caption{Maximal Lyapunov exponent for case \gb\ with different LES-models.}
	\label{tabMaxLyapunovLES}
\end{table}

Table \ref{tabMaxLyapunovLES} lists the MLE for case \gb\ with different LES-models.
The LES-models utilized are the direct filtering approach to the approximate deconvolution model (ADM, \cite{Stolz:1999,Mathew:2003,Foysi:2010a}), based on a $10^{th}$ order low-pass filter, 
a variant of a DMM (Dynamic Mixed Model, see DMM1 in \cite{Martin:2000} for details) 
and the standard and dynamic Smagorinsky model \cite{Martin:2000,Bogey:2006}. 
The constant $C$ present in the Smagorinsky model was either set to fixed values or estimated 
using the dynamic procedure. Using the dynamic procedure the constant reached a maximum of about 
$0.02$, located in the fully turbulent region of the jet in accordance with results from 
\cite{Bogey:2006}.
It can be observed that the MLE varies considerably with the different LES models.
A general tendency is that the MLEs decrease with an increase of the expected dissipation 
of the models. The results suggest, too, that a simulation performed with the DSM can be controlled much more efficiently than a corresponding simulation using the filter variant of the approximate deconvolution model, which is less dissipative. The Smagorinsky models tend to overestimate 
dissipation according to \cite{Martin:2000,Bogey:2006} and thus the increased optimization efficiency is accompanied by a loss in physical accuracy. The results are in line with the observed reduction
of the cost functional in Figure \ref{figJetsInCompRes}, where one could see a better control performance 
of the simulations with less resolution capabilities. They also agree with the control simulations artificially reducing the control space by using only a subset of the control points (Figures \ref{figJetsInCompGapELES3D} and \ref{figJetsInCompGapDNS3D}).

\section{Exact Tangent-Linear Solver by Algorithmic Differentiation}

In order to investigate in detail very large sensitivity values, which have been observed in optimizations with large control horizons, the primal flow solver has been differentiated using Algorithmic Differentiation \citep{GriewankBook} techniques. In this way a tangent-linear solver has been generated, which fully retains the features of the primal flow solver. Moreover, the AD based tangent-linear solver delivers exact sensitivities of the objective function at any convergence level achieved by the primal simulation. Although, the tangent-linear solver requires a run-time that increases linearly with the number of control parameters, it is still a perfect tool to evaluate sensitivities for large control horizons and make an assessment of the situation. For LES/DNS type of problems with a large control horizon, the AD based adjoint solver runs quickly out of memory, and therefore it is not feasible as a sensitivity evaluation method. On the other hand, the tangent-linear method enables sensitivity evaluation of large scale problems at a reasonable memory demand. In the following, we introduce briefly the forward mode of Algorithmic Differentiation, which is used to generate the tangent-linear code used in the study. We also show that the results that are obtained with this code are mathematically equivalent to the adjoint mode. 

Algorithmic Differentiation, also known as Automatic Differentiation, is the name given to a set of techniques developed for the differentiation of functions that are implemented as computer programs. Given a computer code, the AD tool generates a differentiated code that computes the derivatives of the output variables with respect to input variables. The basic idea behind AD is that, the given code can be interpreted as a sequence of simple elementary operations involving additions, multiplications and use of intrinsic functions like $\sin, \exp$ etc. Assuming that the computer code is composed of piece-wise differentiable functions, the derivative code is then obtained by applying the chain rule of differentiation sequentially to each of these elementary operations.

In general, chain rule can be applied in two ways to a given set of elementary operations. The first way, which appears to be more natural, is the so-called forward or tangent-linear mode of AD. Using the forward mode, the chain rule is applied to every operation in a sequence that starts from the input parameters and ends with the output parameters. Therefore, each operation in the data flow is differentiated with respect to a specified direction vector. The resulting derivative expressions are then evaluated simultaneously with the operations of the original function. In contrast to the forward mode, the reverse or adjoint mode of AD applies the chain rule in the reverse order, in which the operations are performed in the original computer program. Note that, both the forward and reverse modes produce exactly the same result.

To simplify the derivation of the tangent-linear method, we consider an objective function over a time interval $[0,T]$ that is to be minimized or maximized
\begin{equation}
J = J(Y_0,Y_1,\ldots,Y_N,X),
\end{equation}
where $Y_i$ is the state vector at the time iteration $i$ and $X$ is the control vector. $N$ is the number of time iterations performed in the primal simulation to reach the final time $t=T$. Note that the above equation is the most general form, in which the control interval starts from $t=0$ and goes up to $t=T$. 
 
Assuming that we have discrete state solutions $Y_0,Y_1,Y_2,\ldots, Y_N$ over the time interval $[0,T]$, the true dependency between the objective function $J$ and the design vector $X$ is given by the relation 

\begin{equation}
J = J(Y_0,Y_1,\ldots,Y_N,X), \quad \mbox{such that} \quad Y_{k+1}=G(Y_k,X),\; k=0,1,\ldots, N-1, 
\label{eq:J}
\end{equation}
where $G$ is a mapping of a state space into itself, i.e., a single time iteration of the flow solver including all the intermediate Runge-Kutta steps of the temporal scheme. In our setting, it includes all the operations within a time iteration, e.g, spatial discretization terms, boundary treatment, filtering etc.   

If we differentiate the objective function $J$ with respect to the design parameter vector $X$, we get
\begin{equation}
\frac{dJ}{dX}= \frac{\partial J}{\partial X} + \frac{\partial J}{\partial Y_0}\frac{d Y_0}{d X} + \frac{\partial J}{ \partial Y_1}\frac{d Y_1}{d X} + \ldots +\frac{\partial J}{\partial Y_N}\frac{d Y_N}{d X}.
\label{eq:dJdX}
\end{equation}
The initial solution $Y_0$ does not depend on the control $X$, so the above equation simplifies to
\begin{equation}
\frac{dJ}{dX}= \frac{\partial J}{\partial X} + \frac{\partial J}{ \partial Y_1}\frac{d Y_1}{d X} + \ldots +\frac{\partial J}{\partial Y_N}\frac{d Y_N}{d X}.
\label{eq:dJ_dX}
\end{equation}

On the other hand, by differentiating the discrete mappings $Y_{k+1}=G(Y_k,X),\; k=1,2,\ldots,N$, we get 

\begin{equation}
\frac{dY_1}{dX} = \frac{\partial G(Y_0,X)}{\partial X}, \quad \frac{dY_{i+1}}{dX} = \frac{\partial G(Y_{i},X)}{\partial Y_{i}} \frac{dY_{i}}{dX} + \frac{\partial G(Y_i,X)}{\partial X},\; i=1,\ldots,N-1.
\label{eq:dYdXchain}
\end{equation}
 
The directional derivation for the given arbitrary differentiation direction $\dot X$ is given by
\begin{equation}
\frac{dJ}{dX} \dot X = \frac{\partial J}{\partial X} \dot X + \sum_{i=1}^N \frac{\partial J}{ \partial Y_i} \frac{d Y_i}{dX} \dot X,
\label{eq:dJ_dX_sum1}
\end{equation}

if we denote the matrix-vector product $\frac{dY_i}{dX} \dot{X}$ by $\dot Y_i$, the above equation can be rewritten as 

\begin{equation}
\frac{dJ}{dX} \dot X = \frac{\partial J}{\partial X} \dot X + \sum_{i=1}^N \frac{\partial J}{ \partial Y_i} \dot Y_i,
\label{eq:dJ_dX_sum2}
\end{equation}
where the $\dot Y_i$ is given by the recursion

\begin{equation}
\dot Y_1 = \frac{\partial G(Y_0,X)}{\partial X} \dot{X}, \quad \dot Y_i = \frac{\partial G(Y_{i},X)}{\partial Y_{i}} \dot Y_{i-1} + \frac{\partial G(Y_i,X)}{\partial X} \dot X,\; i=2,\ldots,N.
\label{eq:Yi_recursion}
\end{equation}

The tangent-linear code that performs the solution procedure given in Eqs. (\ref{eq:dJ_dX_sum2}) and (\ref{eq:Yi_recursion}) can be generated automatically by applying AD techniques on the source code of the primal solver in a black-box fashion. In this way, one obtains a tangent-linear solver that gives exact sensitivity results for a defined forward trajectory of $Y_0,Y_1,\ldots,Y_N$. In the present work, we have used the source transformation tool Tapenade for the differentiation \citep{TapenadeRef13}. The vector $\dot Y_i$ corresponds to the exact linearization of the solution procedure at the time iteration $i$ for the given differentiation direction $\dot X$. One obvious disadvantage of the tangent-linear method, which is outlined above, is the computational cost. Since the forward propagation of derivatives given in Eq. \eqref{eq:dJ_dX_sum2} can be achieved only for a single direction vector $\dot X$ at a time, the procedure must be repeated for all the entries of the gradient vector $dJ/dX$. 

To show that the sensitivities obtained from the tangent-linear solver are equivalent to the adjoint results, one can take the transpose of the Eq. \eqref{eq:dJ_dX} and multiply it with a weight vector $\bar J$ 

\begin{equation}
\left( \frac{dJ}{dX} \right)^\top \!\!\bar J =  \left( \frac{\partial J}{\partial X} \right)^\top \!\!\bar J 
+   \left( \frac{d Y_1}{d X} \right)^\top  \left( \frac{\partial J}{ \partial Y_1} \right)^\top \!\!\bar J  
+ \ldots + \left( \frac{d Y_N}{d X} \right)^\top \left( \frac{\partial J}{ \partial Y_N} \right)^\top\!\! \bar J.
\label{eq:dJdX_transpose}
\end{equation}

Similarly transposing the relations for $dY_i/dX$ yields
\begin{small}
\begin{eqnarray}
\left( \frac{dY_1}{dX} \right)^\top&=& \left( \frac{\partial G(Y_0,X)}{\partial X} \right)^\top, \quad \left( \frac{dY_{i+1}}{dX} \right)^\top =  \left( \frac{dY_{i}}{dX} \right)^\top \frac{\partial G(Y_{i},X)^\top}{\partial Y_{i}}+ \left( \frac{\partial G(Y_i,X)}{\partial X} \right)^\top,\nonumber\\
&&\hfill i=1,\ldots,N-1.
\label{eq:dYdXchain_T}
\end{eqnarray}
\end{small}
Combining the both equations, we get the adjoint sensitivity equation given as
\begin{small}
\begin{eqnarray*}
\left( \frac{dJ}{dX} \right)^{\!\!\top}\!\!\! \bar J\!\!\! & = &\!\!\!  \left( \frac{\partial J}{\partial X} \right)^\top \bar J + \frac{ \partial G(Y_0,X)^\top}{\partial X} \frac{\partial J^\top}{\partial Y_1} \bar J   \nonumber\\
             \!\!\!&+&\!\!\!\left( \frac{\partial G(Y_1,X)^\top}{\partial X} +\frac{\partial G(Y_0,X)^\top}{\partial X} \frac{\partial G(Y_1,X)^\top}{\partial Y_1}  \right) \frac{\partial J^\top}{\partial Y_2} \bar J + \left( \frac{\partial G(Y_2,X)^\top}{\partial X} \right.\nonumber \\
              \!\!\!& + &\!\!\!   \left. \frac{\partial G(Y_1,X)^\top}{\partial X} \frac{\partial G(Y_2,X)^\top}{\partial Y_2} + \frac{\partial G(Y_0,X)^\top}{\partial X}  \frac{\partial G(Y_1,X)^\top}{\partial Y_1}  \frac{\partial G(Y_2,X)^\top}{\partial Y_2} \right) \!\!\frac{\partial J^\top}{\partial Y_3}\! \bar J \nonumber \\
              \!\!\!& + &\!\!\! \ldots \nonumber \\
              \!\!\!& + &\!\!\! \left( \frac{\partial G(Y_{N-1},X)^{\!\!\top}}{\partial X} + \frac{\partial G(Y_{N-2},X)^{\!\!\top}}{\partial X} \frac{\partial G(Y_{N-1},X)^{\!\!\top}}{\partial Y_{N-1}} +  \ldots \right. \nonumber \\
              \!\!\!&+  &\!\!\! \frac{\partial G(Y_1,X)^{\!\!\top}}{\partial X} \frac{\partial G(Y_2,X)^{\!\!\top}}{\partial Y_2}  \!\ldots\! \frac{\partial G(Y_{N-2},X)^{\!\!\top}}{\partial Y_{N-2}} \frac{\partial G(Y_{N-1},X)^{\!\!\top}}{\partial Y_{N-1}}  \nonumber\\
              \!\!\!&+&\!\!\!\left.  \frac{\partial G(Y_0,X)^{\!\!\top}}{\partial X}\frac{\partial G(Y_1,X)^{\!\!\top}}{\partial Y_1} \!\ldots\! \frac{\partial G(Y_{N-3},X)^{\!\!\top}}{\partial Y_{N-3}} \frac{\partial G(Y_{N-2},X)^{\!\!\top}}{\partial Y_{N-2}} \frac{\partial G(Y_{N-1},X)^{\!\!\top}}{\partial Y_{N-1}}     \right)\!\! \frac{\partial J^{\!\!\top}}{\partial Y_N} \bar J.
\label{eq:dJ_dX_long_T}
\end{eqnarray*}
\end{small}
Similar to the tangent-linear solver, the adjoint code that performs the solution procedure given in the above equation can be generated automatically by applying AD techniques on the source code of the primal solver in a black-box fashion, with the only difference that, this time the differentiation must be done in reverse mode. Since the objective function is a scalar, the weight vector $\bar J$ is also a scalar and can be simply chosen as $1$. In this way, the complete gradient vector $dJ/dX$ can be evaluated only with a single run of the adjoint code. The memory demand, on the other hand, increases linearly with the number of time iterations performed in the primal simulation as the state vector $Y$ must be available in the adjoint evaluation in the reverse order, i.e, $Y_{N-1},Y_{N-2}, \ldots, Y_1$.   

Note that from the above equations, we can easily derive the relationship 
\begin{equation}
\bar J \dot X = \dot J \bar X,
\end{equation}
which shows the relationship between the adjoint control sensitivities $\bar X$ and the directional derivative $\dot J$ in the direction $\dot X$. In conclusion, by applying AD techniques it is guaranteed that the tangent-linear sensitivity results are equivalent to the adjoint results.

\section{Tangent-Linear Sensitivities with Long Time Horizon}

For the validation of the AD based tangent-linear solver, the round jet configuration with $Re=10000$ and $Ma=0.9$ has been used. As far as the computational grid is concerned, a structured grid with $448 \times 216 \times 288$ grid points are used (JetG2 test case). For the forcing term, to simplify the analysis, we take a scalar term that is added to the energy equation
\begin{equation}
\frac{\partial p}{\partial t}  =  -\frac{\partial}{\partial x_i} p u_i + \frac{\partial}{\partial x_i} \lambda (\gamma - 1) \frac{\partial}{\partial x_i} - (\gamma - 1 ) p \frac{\partial}{\partial x_i} u_i + (\gamma - 1 ) \tau_{i j} \frac{\partial}{\partial x_j} u_i + \rho R s(x,t) g,
\end{equation}
where $g$ is the scalar forcing term and $s$ is the windowing function to ensure a smooth transition from uncontrolled to
controlled areas in the flow domain. The windowing function is given by 
\begin{equation}
s(x,t) = s_{window}(x,2\Delta_x)s_{window}(z,2\Delta_z)s_{window}(t,5\Delta_t)
\end{equation}
and
\begin{equation}
s_{window}(k,\Delta_x) = \frac{1}{2} ( erf((k-k_{start}-2\Delta)/\Delta) - erf((k-k_{end}+2\Delta)/\Delta) ),
\end{equation}
where $erf(x)$ is the error function, $\Delta$ is the grid spacing for spatial and the time-step for temporal directions and $k_{start}$ and $k_{end}$ are the start and end positions/times of the
controlled area. With this definitions, the control can be interpreted as a temperature forcing term added to the energy equation.

\begin{figure}
\begin{minipage}{0.5\textwidth}
\includegraphics[scale=0.4]{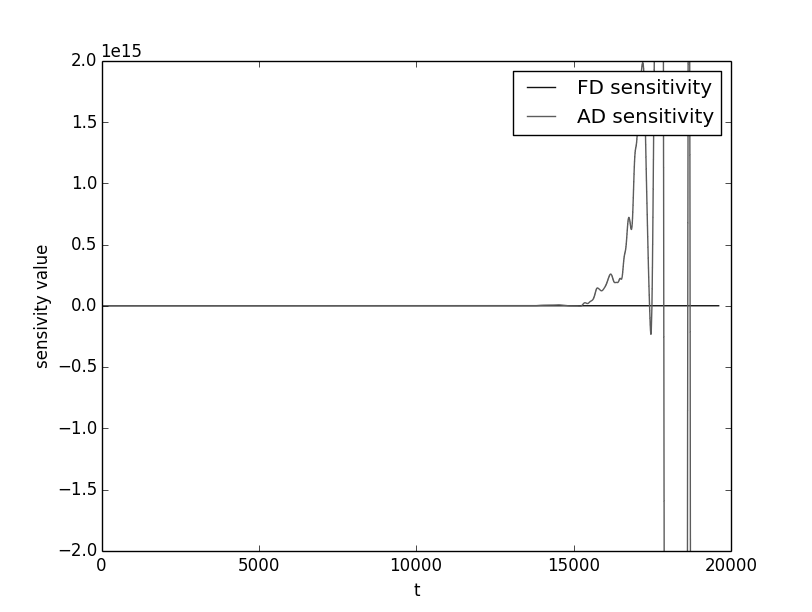} 
\end{minipage}
\begin{minipage}{0.5\textwidth}
\includegraphics[scale=0.4]{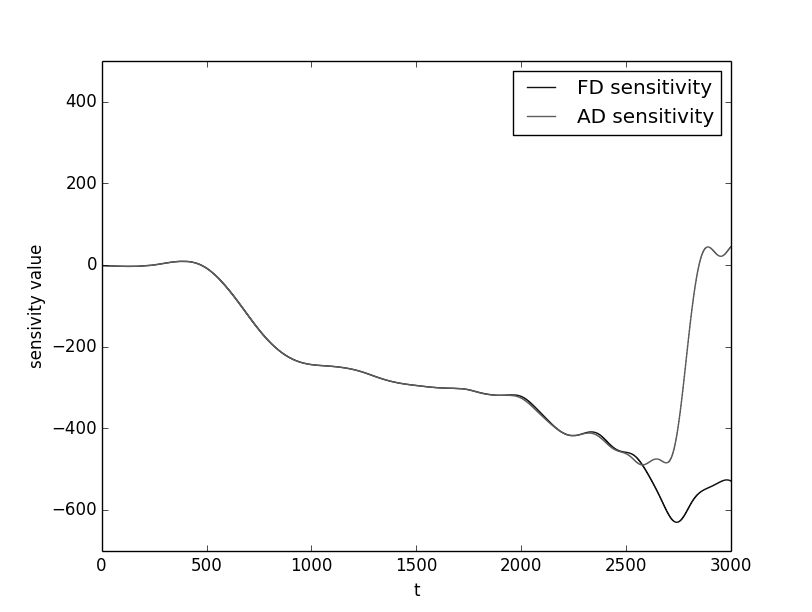} 
\end{minipage}
\caption{Comparison of AD and FD sensitivities (Left: complete time horizon, right: initial $3000$ iterations)}
\label{fig:ADvsFD}
\end{figure}

In Figure \ref{fig:ADvsFD}, the time histories of the tangent-linear sensitivities obtained by AD and forward finite differences are shown. It can be observed that, in contrast to the AD sensitivities, finite difference results remain bounded and do not tend to overflow. In the right figure, the same trajectories up to $3000$ iterations are shown. From the figure, it can observed that both curves coincide each other perfectly within this narrow time window, which is an indication that the tangent-linear sensitivities are correct. It is interesting that, at around $2500$th time iteration, the tangent-linear results start to deviate significantly from the finite difference results. In other words, two trajectories bifurcate. At the later iterations, the tangent-linear sensitivity values simply grow because of the "butterfly effect". The FD sensitivities, however, do not show this behavior and they are, therefore, certainly unreliable after a certain number of time iterations.

%
%
%

We can conclude that the sensitivity trajectory obtained from the LES simulation shows a very similar behavior to the Lorenz system \citep{Lorenz}. In both cases, the perturbed trajectory of the objective function stays very close to the original trajectory in the initial time steps. The gap between the both trajectories, however, tend to increase as more time iterations are taken in the forward-in-time integration. In both cases, the sensitivity results suffer from the so-called "butterfly" effect and tend to overflow. The objective function values, on the other hand, remain bounded. Unfortunately, this phenomena renders gradient based approaches for flow control problems almost useless as the gradient information becomes unreliable if the control horizon is taken large. On the other hand, choosing a small time horizon may not be appropriate as well as the resolution of the physics is concerned. Usage of AD techniques, unfortunately, does not alleviate this problem. On the contrary, it may even worsen the situation since all kind of instabilities present in the solution (either physical or numerical) are exactly differentiated by the underlying AD library. Unfortunately, differentiating noisy parts of the primal solution corrupts the sensitivity results as the exact derivative values tend to go to infinity. Using more advanced regularization techniques like Least Square Shadowing (LSS) \citep{WANG2014} is a promising method to overcome this problem. It comes, however, at the expense of increased computational cost.

\section{Conclusions}
In this paper, results of several adjoint based flow control studies performed on turbulent plane and round jet configurations are presented. For the flow simulations, a LES/DNS finite difference solver with a high order spatial scheme has been used. Based on this flow solver, continuous and discrete adjoint versions have been developed to able to evaluate gradient vectors efficiently. From the results, it has been observed that the cost functional can be significantly reduced using an adjoint based L-BFGS optimization algorithm provided that the control horizons are kept small enough. A significant difference between continuous and discrete approaches has not been observed as both approaches are good enough to evaluate accurate sensitivity gradients. For large control horizons, on the other hand, it has been observed that flow control is no more possible as the sensitivity gradients tend to grow in time. To be able to make a better assessment of the situation and exclude the effect of all possible consistency errors in the adjoints, the same flow solver has been differentiated using the machine accurate AD techniques in forward mode. In this way, an exact tangent-linear solver has been developed that corresponds to an exact linearization of the flow solver with all its underlying features. The results that are obtained with the tangent-linear solver gave a similar picture. In the initial time iterations, the tangent-linear results matched perfectly with the finite difference approximations. After a certain time iteration, however, the both trajectories started to deviate from each other. Similar to the adjoint results, tangent-linear sensitivities also grow rapidly with the increasing simulation time. From the results, it can be concluded that after a certain size of the control horizon adjoint based flow control is no more viable without regularization techniques.  





\section*{References}
\begin{thebibliography}{99}


\bibitem{nielsen_2011} E.J. Nielsen, W.T. Jones, ``Integrated design of an active flow control system using a time-dependent adjoint method", \textit{Mathematical Modeling of Natural Phenomena}, Vol. 6, Number 3, (2011), pp. 141-165.

\bibitem{Bewley:2001} T. Bewley, P. Moin, R. Temam, ``DNS-based predictive control of turbulence: an optimal benchmark for feedback algorithms", \textit{Journal Fluid Mechanics}, Vol. 447, (2001), pp. 179-225. doi:https://doi.org/10.1017/S0022112001005821

\bibitem{Wei:2006} M. Wei, J.B Freund, ``A noise-controlled free shear flow", \textit{Journal of Fluid Mechanics}, Vol. 546, (2006), pp. 123-152. doi:https://doi.org/10.1017/S0022112005007093

\bibitem{Jameson:1998} A. Jameson, L. Martinelli, N.A. Pierce, ``Optimum aerodynamic design using the {N}avier-{S}tokes equations", \textit{Theoretical and Computational Fluid Dynamics}, Vol. 10,  Issue 1-4, (1998), pp. 213-237. doi:https://doi.org/10.1007/s001620050060


\bibitem{Nielsen} E.J. Nielsen, B. Diskin, N.K. Yamaleev, ``Discrete adjoint-based design optimization of unsteady turbulent flows on dynamic unstructured grids", \textit{AIAA Journal}, Vol. 48(6), (2010), pp. 1195-1206. doi:https://doi.org/10.2514/1.J050035


\bibitem{Pet2010} J.E.V. Peter, R.P. Dwight, ``Numerical sensitivity analysis for aerodynamic optimization: A survey of approaches and applications", \textit{Computers \& Fluids}, Vol. 39(3), (2010), pp. 373–391. doi:https://doi.org/10.1016/j.compfluid.2009.09.013

\bibitem{anil2013} A. Nemili, E. \"{O}zkaya, N.R. Gauger, F. Kramer, A. Carnarius, F. Thiele, ``Discrete adjoint based sensitivity analysis for optimal flow control of a $3${D} high-lift configuration", \textit{AIAA Paper 2013-2585}, (2013). doi: https://doi.org/10.2514/6.2013-2585 

\bibitem{HascoetGK09} L. Hasco\"et, ``Reversal Strategies for Adjoint Algorithms", \textit{From Semantics to Computer Science. Essays in memory of {G}illes {K}ahn}, Cambridge University Press, (2009). pp. 487-503. 


\bibitem{Mani08} K. Mani, D.J. Mavriplis, ``Unsteady discrete adjoint formulation for two-dimensional flow problems with deforming meshes", \textit{AIAA Journal}, Vol. 46(6), (2008), pp. 1351-1364, doi:https://doi.org/10.2514/1.29924 


\bibitem{Nada02} S.K. Nadarajah, A. Jameson, ``Optimal control of unsteady flows using a time accurate method", \textit{AIAA Paper 2002-5436}, (2002). doi: https://doi.org/10.2514/6.2002-5436  

\bibitem{Rump07} M.P. Rumpfkeil, D.W. Zingg, ``A general framework for the optimal control of unsteady flows with applications", \textit{AIAA Paper 2007-1128}, 2007. doi: https://doi.org/10.2514/6.2007-1128  


\bibitem{Griewank2000} A. Griewank, A. Walther, ``Algorithm 799: Revolve: An implementation of checkpointing for the reverse or adjoint mode of computational differentiation", \textit{ACM Trans. Math. Softw.}, Vol. 26(1), (2000), pp. 19-45. doi:http://doi.acm.org/10.1145/347837.347846

\bibitem{Stumm09} P. Stumm, A. Walther, ``Multistage approaches for optimal offline checkpointing", \textit{SIAM Journal on Scientific Computing}, Vol. 31(3), (2009), pp. 1946-1967. doi:http://doi.acm.org/10.1137/080718036


\bibitem{Zhou_etal_2016a} B.Y. Zhou, N.R. Gauger, S.R. Koh, M. Meinke, W. Schr\"{o}der, ``A Discrete adjoint framework for trailing-edge turbulence control and noise minimization via porous material", \textit{AIAA Paper 2016-2777}, (2016). doi: https://doi.org/10.2514/6.2007-1128  

\bibitem{NeOeGaKrTh2016b} A. Nemili, E. \"{O}zkaya, N.R. Gauger, F. Kramer, F. Thiele, ``Aerodynamic Optimisation of active flow control on a three-element high-lift configuration", \textit{Proceedings of 9th International Conference on Computational Fluid Dynamics}, ICCFD9-2016-160, (2016). available online at http://iccfd9.itu.edu.tr/assets/pdf/papers/ICCFD9-2016-160.pdf 


\bibitem{wang13} Q. Wang, ``Forward and adjoint sensitivity computation of chaotic dynamical systems", \textit{Journal of Computational Physics}, Vol. 235, (2013), pp. 1-13. doi:http://dx.doi.org/10.1016/j.jcp.2012.09.007

\bibitem{Ozk16Arxiv} E. \"{O}zkaya, N.R. Gauger, A. Nemili, ``Chaotic behavior of stiff {ODE}s and their derivatives: An illustrative example", \textit{arXiv preprint 1610.03358}, 

\bibitem{WANG2014} Q. Wang, R. Hu, P. Blonigan, ``Least Squares Shadowing sensitivity analysis of chaotic limit cycle oscillations", \textit{Journal of Computational Physics}, Vol. 267, (2014), pp. 210-224. doi:https://doi.org/10.1016/j.jcp.2014.03.002

\bibitem{SteffiOS} S. G\"{u}nther, N.R. Gauger, Q. Wang, ``A framework for simultaneous aerodynamic design optimization in the presence of chaos", \textit{Journal of Computational Physics}, Vol. 328, (2017), pp. 387-398. https://doi.org/10.1016/j.jcp.2016.10.043

\bibitem{Guenther2017} S. G\"{u}nther, N.R. Gauger, J. Schroder, ``A non-intrusive parallel-in-time adjoint solver with the {XBraid} library", \textit{arXiv preprint 1705.00663}, https://arxiv.org/pdf/1705.00663.pdf

\bibitem{Marinc2012} D. Marinc, H. Foysi, ``Investigation of a continuous adjoint-based optimization procedure for aeroacoustic control of plane jets", \textit{International Journal of Heat and Fluid Flow}, Vol. 38, (2012), pp. 200-212. doi: 
https://doi.org/10.1016/j.ijheatfluidflow.2012.07.005

\bibitem{Hu:1996}  F.Q. Hu, M.Y. Hussaini, J.L. Manthey, ``Low-dissipation and lowdispersion Runge-Kutta schemes for computational acoustics", \textit{J. Comput. Phys.}, Vol. 124, (1996), pp. 177-197. doi: https://doi.org/10.1006/jcph.1996.0052

\bibitem{Lodato2008} G. Lodato, P. Domingo, L. Vervisch, ``Three-dimensional boundary conditions for direct and large-eddy simulation of compressible viscous flows", \textit{Journal of Computational Physics}, Vol. 227, (2008), pp. 5105 - 5143. doi: https://doi.org/10.1016/j.jcp.2008.01.038

\bibitem{Foysi:2010a} H. Foysi, M. Mellado, S. Sarkar, ``Simulation and comparison of variable density round and plane jets", \textit{International Journal of Heat \& Fluid Flow}, Vol. 31, (2010), doi: https://doi.org/10.1016/j.ijheatfluidflow.2009.12.001

\bibitem{Stolz:1999} S. Stolz, N.A. Adams, ``An approximate deconvolution procedure for large-eddy simulation", \textit{Phys. Fluids}, Vol. 11, Issue 7, (1999), pp. 1699-1701, doi: https://doi.org/10.1063/1.869867

\bibitem{Mathew:2003} J. Mathew, R. Lechner, H. Foysi, J. Sesterhenn, R. Friedrich, ``An explicit filtering method for large eddy simulation of compressible flows", \textit{Phys. Fluids}, Vol. 15(8), (2003), pp. 2279-2289, doi: https://doi.org/10.1063/1.1586271

\bibitem{Mathew:2006} J. Mathew, H. Foysi, R. Friedrich, ``A new approach to LES based on explicit filtering", \textit{Int. J. Heat Fluid Flow}, Vol. 27(4), (2006), pp. 594-602, doi: https://doi.org/10.1016/j.ijheatfluidflow.2006.02.007


\bibitem{Petersen} M.R. Petersen, D. Livescu, ``Forcing for statistically stationary compressible isotropic turbulence", \textit{Phys. of Fluids}, Vol. 22(11), (2010), pp. 307-314. doi: https://doi.org/10.1063/1.3488793


\bibitem{Kuptsov} P. Kuptsov, U. Parlitz, ``Theory and computation of covariant Lyapunov vectors", \textit{J. of Nonl. Science}, Vol. 22(5), (2012), pp. 727-762. doi: 
https://doi.org/10.1007/s00332-012-9126-5

\bibitem{Kim:2010} J. Kim, D.J. Bodony, J.B. Freund, ``A high-order, overset-mesh algorithm for adjoint-based optimization for aeroacoustics control", \textit{AIAA 2010-3818}, (2010). doi: https://doi.org/10.2514/6.2010-3818 

\bibitem{Pando} M.F. de Pando, D. Sipp, P.J. Schmid, ``Efficient evaluation of the direct and adjoint linearized dynamics from compressible flow solvers", \textit{J. Comput. Phys.}, Vol. 231(23), (2012), pp. 7739-7755. doi:https://doi.org/10.1016/j.jcp.2012.06.038


\bibitem{Collis:2000} S.S. Collis, Y. Chang, S. Kellog, R.D. Prabhu, ``Large eddy simulation and turbulence control", \textit{AIAA Paper 2000-2564}, (2000). 


\bibitem{Marschner:1994} S.R. Marschner, R.J Lobb, ``An evaluation of reconstruction filters for volume rendering", \textit{Proceedings of the conference on Visualization'94}, (1994), pp. 100-107. 

\bibitem{Pope} S.B. Pope, ``Turbulent Flows", \textit{Cambridge University Press}, 2000.

\bibitem{Bogey:2006} C. Bogey, C. Bailly, ``Large eddy simulation of round free jets using explicit filtering with/without dynamic {S}magorinsky model", \textit{Int. J. Heat and Fluid Flow}, Vol. 27(4), (2006), pp. 603-610. doi:https://doi.org/10.1016/j.ijheatfluidflow.2006.02.008

\bibitem{Crisanti:1993} Crisanti, A., Jensen, M.H., Vulpiani, A., Paladin, G., ``Intermittency and predictability in turbulence", \textit{Phys. Rev. Lett.}, Vol. 70(2), (1993), pp. 166-169. doi: 10.1103/PhysRevLett.70.166

\bibitem{Martin:2000} M.P. Martin, U. Piomelli, G.V.Candler, ``Subgrid-scale models for compressible large-eddy simulations", \textit{Theoret. Comp. Fluid Dyn.}, Vol. 13, (2000), pp. 361-376. doi: https://doi.org/10.1007/PL00020896

\bibitem{GriewankBook} A.Griewank, A. Walther, Evaluating Derivatives: Principles and Techniques of Algorithmic Differentiation, second ed. Society for Industrial and Applied Mathematics, 2000.

\bibitem{TapenadeRef13} L. Hasco{\"e}t, V. Pascual, ``The {T}apenade {A}utomatic {D}ifferentiation tool: {P}rinciples, {M}odel, and {S}pecification", \textit{{ACM} {T}ransactions {O}n {M}athematical {S}oftware}, Vol. 39(3), (2013), pp. 20:1-20:43. doi: http://dx.doi.org/10.1145/2450153.2450158 

\bibitem{Lorenz} E.N. Lorenz, ``Deterministic nonperiodic flow ", \textit{Journal of the Atmospheric Sciences}, Vol. 20(2), (1963), pp. 130-141. 






 \end{thebibliography}


\section*{Acknowledgments}
This research was funded by the German Research Foundation (DFG) under the project numbers, $FO \; 674/4-1$ and $GA \; 857/9-1$. The authors gratefully acknowledge the computing time granted by the Allianz f\"{u}r Hochleistungsrechnen Rheinland-Pfalz (AHRP), Germany.

\end{document}